\documentclass[11pt]{amsart}

\usepackage{amsmath}

\usepackage{amssymb}

\newcommand{\C}{{\mathbb C}}       
\newcommand{\R}{{\mathbb R}}       
\newcommand{\ds}{\displaystyle }

\newcommand{\rf}[1]{{(\ref{#1})}}

\newcommand{\supp}{{\rm supp}}

\newcommand{\vphi}{{\varphi}}
\newcommand{\ve}{{\varepsilon}}

\newcommand{\bmo}{{B\!M\!O}}

\newtheorem{theorem}{Theorem}[section]
\newtheorem{lemma}[theorem]{Lemma}

\theoremstyle{definition}

\newtheorem{example}[theorem]{Example}

\theoremstyle{remark}
\newtheorem{rem}[theorem]{Remark}

\numberwithin{equation}{section}

\newcommand{\brem}{\begin{rem}}
\newcommand{\erem}{\end{rem}}

\newcommand{\bexam}{\begin{example}}
\newcommand{\eexam}{\end{example}}

\begin{document}

\title[Weak (1,1) inequality with non doubling 
measures]{A proof of the weak (1,1) inequality for singular integrals 
with non doubling measures based on a Calder\'on-Zygmund
decomposition}

\author[XAVIER TOLSA]{Xavier Tolsa}

\address{Department of Mathematics, Chalmers, 412 96 G\"oteborg, Sweden}

\email{xavier@math.chalmers.se}

\thanks{Supported by a postdoctoral grant from the European Comission for the
TMR Network ``Harmonic Analysis''. Also
partially supported by grants DGICYT PB96-1183 and CIRIT
1998-SGR00052 (Spain).}


\date{February 4, 2000.}


\begin{abstract}
Given a doubling measure $\mu$ on $\R^d$, it is a
classical result of harmonic analysis that Calder\'on-Zygmund operators
which are bounded in $L^2(\mu)$ are also of weak type $(1,1)$.
Recently it has been shown that the same result holds if one substitutes 
the doubling condition on $\mu$ by a mild growth condition on $\mu$.
In this paper another proof of this result is given. The proof 
is very close in spirit to the classical argument for doubling measures and it
is based on a new Calder\'on-Zygmund decomposition adapted to the non doubling 
situation.
\end{abstract}

\maketitle


\section{Introduction} \label{sec1}

Let $\mu$ be a positive Radon measure on $\R^d$ satisfying the growth condition
\begin{equation}  \label{creix}
\mu(B(x,r))\leq C_0\,r^n\quad \mbox{for all $x\in \R^d,\, r>0,$}
\end{equation}
where $n$ is some fixed number with $0<n\leq d$.
We do {\em not} assume that $\mu$ is doubling
[$\mu$ is said to be doubling if there exists some constant $C$
such that $\mu(B(x,2r)) \leq C\,\mu(B(x,r))$ for all $x\in \supp(\mu)$, $r>0$].
Let us remark that the doubling condition on the underlying measure $\mu$ 
on $\R^d$ is an essential assumption in most
results of classical Calder\'on-Zygmund theory. However, recently it
has been shown that a big part of the classical theory remains valid
if the doubling assumption on $\mu$ is substituted by the size condition
\rf{creix} (see for example the references cited at the end of the paper).

In this note we will prove that Calder\'on-Zygmund operators (CZO's) 
which are bounded
in $L^2(\mu)$ are also of weak type $(1,1)$, as in the usual doubling
situation. This result
has already been proved in \cite{Tolsa1} in the particular case of the Cauchy
integral operator, and by Nazarov, Treil and Volberg \cite{NTV2} in the general
case. The proof that we will show here is different from the one of \cite{NTV2}
(and also from the one of \cite{Tolsa1}, of course) and it is closer
in spirit to the classical proof of the corresponding result for doubling 
measures. The basic tool for the proof is a decomposition 
of Calder\'on-Zygmund type for functions in $L^1(\mu)$ obtained in 
\cite{Tolsa4}.

Our purpose in writing this paper is not only to obtain another proof in the non
doubling situation of the basic result that CZO's bounded in $L^2(\mu)$ are of
weak type $(1,1)$, but to show that the Calder\'on-Zygmund decompositon of
\cite{Tolsa4} is a good substitute of its classical doubling version.

Let us introduce some notation and definitions. A kernel
$k(\cdot,\cdot)$ from 
$L^1_{loc}((\R^d\times\R^d) \setminus \{(x,y) :x=y\})$ is called a
Calder\'on-Zygmund kernel if
\begin{enumerate}
\item $\ds |k(x,y)|\leq \frac{C}{|x-y|^n},$
\item there exists $0<\delta\leq1$ such that
$$|k(x,y)-k(x',y)| + |k(y,x)-k(y,x')|
\leq C\,\frac{|x-x'|^\delta}{|x-y|^{n+\delta}}$$
if $|x-x'|\leq |x-y|/2$.
\end{enumerate}
Throughout all the paper we will assume that $\mu$ is a Radon measure on
$\R^d$ satisfying \rf{creix}. The CZO associated to the kernel
$k(\cdot,\cdot)$ and the measure $\mu$ is defined (at least, formally) as
$$Tf(x) = \int k(x,y)\,f(y)\,d\mu(y).$$
The above integral may be not convergent for many functions $f$ because
$k(x,y)$ may have a singularity for $x=y$. For this reason,
one introduces the truncated operators $T_\varepsilon$, $\varepsilon>0$:
$$T_\varepsilon f(x) = \int_{|x-y|>\varepsilon} k(x,y)\,f(y)\,d\mu(y),$$
and then one says that $T$ is bounded in $L^p(\mu)$ if the operators
$T_\varepsilon$ are bounded in $L^p(\mu)$ uniformly on $\varepsilon>0$. 
It is said that $T$ is bounded from $L^1(\mu)$ into $L^{1,\infty}(\mu)$
(or of weak type $(1,1)$) if
$$\mu\{x:|T_\ve f(x)|>\lambda\} \leq C\, \frac{\|f\|_{L^1(\mu)}}{\lambda}$$
for all $f\in L^1(\mu)$, uniformly on $\ve>0$. Also, $T$ is bounded
from $M(\C)$ (the space of complex Radon measures) into $L^{1,\infty}(\mu)$ if
$$\mu\{x:|T_\ve\nu(x)|>\lambda\} \leq C\, \frac{\|\nu\|}{\lambda}$$
for all $\nu\in M(\C)$, uniformly on $\ve>0$. In the last inequality,
$T_\ve\nu(x)$ stands for $\int_{|x-y|>\varepsilon} k(x,y)\,d\nu(y)$ and
$\|\nu\|\equiv |\nu|(\R^d)$.

The result that we will prove in this note is the following.

\begin{theorem}\label{teo}
Let $\mu$ be a Radon measure on $\R^d$ satisfying the growth condition
\rf{creix}. If $T$ is a Calder\'on-Zygmund operator which is bounded in
$L^2(\mu)$, then it is also bounded from $M(\C)$ into
$L^{1,\infty}(\mu)$. In particular, it is of weak type $(1,1)$.
\end{theorem}


\section{The proof}

First we will introduce some additional notation and terminology. 
As usual, the letter $C$ will denote a
constant which may change its value from one occurrence to another. 
Constants with subscripts, such as $C_0$, do not change in different
occurrences.

By a cube $Q\subset R^d$ we
mean a closed cube with sides parallel to the axes. We denote its side length by
$\ell(Q)$ and its center by $x_Q$. 
Given $\alpha>1$ and $\beta>\alpha^n$, we say that $Q$
is $(\alpha,\beta)$-doubling if $\mu(\alpha Q) \leq \beta\,\mu(Q)$, where
$\alpha Q$ is the cube concentric with $Q$ with side length $\alpha\,\ell(Q)$.
For definiteness, if $\alpha$ and $\beta$ are not specified, by a doubling cube
we mean a $(2,2^{d+1})$-doubling cube.

Before proving Theorem \ref{teo} we state some remarks about the existence of
doubling cubes.

\begin{rem} \label{rem1}
Because $\mu$ satisfies the growth condition \rf{creix}, there
are a lot ``big'' doubling cubes. To be precise, given any point
$x\in\supp(\mu)$ and $c>0$, there exists some
$(\alpha,\beta)$-doubling cube $Q$ centered at $x$ with $l(Q)\geq c$. This
follows easily from \rf{creix} and the fact that $\beta>\alpha^n$.
Indeed, if there are no doubling cubes centered at $x$ with $l(Q)\geq c$,
then $\mu(\alpha^n Q)> \beta^n \mu(Q)$ for each $n$, and letting $n\to\infty$
one sees that \rf{creix} cannot hold.
\end{rem}

\begin{rem}\label{rem2}
There are a lot of ``small'' doubling cubes too: if 
$\beta>\alpha^d$, then for $\mu$-a.e. $x\in\R^d$ there
exists a sequence of $(\alpha,\beta)$-doubling cubes $\{Q_k\}_k$ centered
at $x$ with $\ell(Q_k)\to0$ as $k\to \infty$. This is a property that any
Radon measure on $\R^d$ satisfies (the growth condition \rf{creix} is not
necessary in this argument). The proof is an easy exercise on geometric measure
theory that is left for the reader.

Observe that, by the Lebesgue differentiation theorem, for $\mu$-almost all
$x\in\R^d$ one can find
a sequence of $(2,2^{d+1})$-doubling cubes $\{Q_k\}_k$
centered at $x$ with $\ell(Q_k)\to0$ such that
$$\lim_{k\to\infty} \frac{1}{\mu(Q_k)} \int_{Q_k} f \,d\mu = f(x).$$
As a consequence, for any fixed $\lambda>0$, for $\mu$-almost all $x\in\R^d$ such that
$|f(x)|>\lambda$, there exists a sequence of cubes
$\{Q_k\}_k$ centered at $x$ with $\ell(Q_k)\to0$ such that
$$\limsup_{k\to\infty} \frac{1}{\mu(2Q_k)} \int_{Q_k} |f| \,d\mu >
\frac{\lambda}{2^{d+1}}.$$
\end{rem}

In the following lemma we will prove an easy but essential estimate which 
will be used below.
This result has already appeared in previous works (\cite{DM}, \cite{NTV2})
and it plays a basic role in \cite{Tolsa2} and \cite{Tolsa4} too.

\begin{lemma}\label{basic}
If $Q\subset R$ are concentric cubes such that there are no 
$(\alpha,\beta)$-doubling cubes (with $\beta>\alpha^n$) of
the form $\alpha^kQ$, $k\geq0$, with $Q\subset \alpha^kQ\subset R$, then,
$$\int_{R\setminus Q} \frac{1}{|x-x_Q|^n}\,d\mu(x) \leq C_1,$$
where $C_1$ depends only on $\alpha, \beta$, $n$, $d$ and $C_0$.
\end{lemma}

\begin{proof}
Let $N$ be the least integer such that $R\subset \alpha^N Q$. For $0\leq k\leq
 N$ we have $\mu(\alpha^k Q) \leq \mu(\alpha^N Q) /\beta^{N-k}$. Then,
\begin{eqnarray*}
\int_{R\setminus Q} \frac{1}{|x-x_Q|^n}\,d\mu(x) & \leq &
\sum_{k=1}^N \int_{\alpha^k Q \setminus \alpha^{k-1}Q}
\frac{1}{|x-x_Q|^n}\,d\mu(x)\\
& \leq & C\,\sum_{k=1}^N \frac{\mu(\alpha^k Q)}{\ell(\alpha^k Q)^n} \\
& \leq & C\,\sum_{k=1}^{N} 
\frac{\beta^{k-N}\,\mu(\alpha^N Q)}{\alpha^{(k-N)n}\,\ell(\alpha^N Q)^n} \\
& \leq & C\,\frac{\mu(\alpha^N Q)}{\ell(\alpha^N Q)^n} \sum_{j=0}^{\infty}   
\left(\frac{\alpha^n}{\beta}\right)^j \,\leq \, C.
\end{eqnarray*}
\end{proof}

The Calder\'on-Zygmund decomposition mentioned above has been obtained
in Lemma 7.3 of \cite{Tolsa4} and in that paper it has been used to show that 
if a linear operator
is bounded from a suitable space of type $H^1$ into $L^1(\mu)$ and from
$L^\infty(\mu)$ into a space of type $BMO$, then it is bounded in $L^p(\mu)$.
We will use a slight variant of this decompositon to prove Theorem \ref{teo}.
Let us state the result that we need in detail.

\begin{lemma}[Calder\'on-Zygmund decomposition]  \label{lem}
Assume that $\mu$ satisfies \rf{creix}. For any $f\in L^1(\mu)$ and any 
$\lambda>0$ (with $\lambda>2^{d+1}\,\|f\|_{L^1(\mu)}/\|\mu\|$ if 
$\|\mu\|<\infty$) we have:
\begin{itemize}
\item[(a)] There exists a finite family of
almost disjoint cubes $\{Q_i\}_i$ such that
\begin{equation}  \label{cc1}
\frac{1}{\mu(2Q_i)} \int_{Q_i} |f|\,d\mu >\frac{\lambda}{2^{d+1}},
\end{equation}
\begin{equation}  \label{cc2}
\frac{1}{\mu(2\eta Q_i)} \int_{\eta Q_i} |f|\,d\mu \leq
\frac{\lambda}{2^{d+1}} \quad
\mbox{for $\eta >2$,}
\end{equation}
\begin{equation}  \label{cc3}
|f|\leq \lambda \quad \mbox{a.e. ($\mu$) on $\R^d\setminus\bigcup_i
Q_i$}.
\end{equation}

\item[(b)] For each $i$, let $R_i$ be a $(6,6^{n+1})$-doubling cube concentric
with $Q_i$, with $l(R_i)>4l(Q_i)$ and denote
$w_i= \frac{\chi_{Q_i}}{\sum_k \chi_{Q_k}}$. Then,
there exists a family of functions
$\vphi_i$ with $\supp(\vphi_i)\subset R_i$ and with constant sign satisfying
\begin{equation}  \label{cc4}
\int \vphi_i \,d\mu = \int_{Q_i} f\,w_i\,d\mu,
\end{equation}
\begin{equation}  \label{cc5}
\sum_i |\vphi_i| \leq B\,\lambda
\end{equation}
(where $B$ is some constant), and
\begin{equation}  \label{cc6}
\|\vphi_i\|_{L^\infty(\mu)} \,\mu(R_i)\leq
C\, \int_{Q_i}|f|\,d\mu.
\end{equation}

\end{itemize}
\end{lemma}

Let us remark that other related decompositons with non doubling measures have
been obtained in \cite{NTV2} and \cite{MMNO}. However, these results are not
suitable for our purposes.

Although the proof of the lemma can be found in \cite{Tolsa4},
for the reader's convenience we have included it in the last section of the 
present paper.

\begin{proof}[\bf Proof of Theorem \ref{teo}]
We will show that $T$ is of weak type $(1,1)$. By similar arguments, one gets
that $T$ is bounded from $M(\C)$ into $L^{1,\infty}(\mu)$. In this case, one has
to use a version of the Calder\'on-Zygmund decomposition in the lemma above
suitable for complex measures. 

For simplicity we assume $\|\mu\|=\infty$.
Let $f\in L^1(\mu)$ and $\lambda>0$. Let $\{Q_i\}_i$ be the almost disjoint 
family of cubes of Lemma \ref{lem}. 
Let $R_i$ be the smallest $(6,6^{n+1})$-doubling cube of the form $6^kQ_i$,
$k\geq1$. Then we can write $f=g+b$, with
$$g= f\,\chi_{\R^d\setminus\bigcup_i Q_i} + \sum_i \vphi_i$$
and
$$b = \sum_i b_i := \sum_i \left(w_i\, f - \vphi_i\right),$$
where the functions $\vphi_i$ satisfy \rf{cc4}, \rf{cc5} \rf{cc6} and 
$w_i= \frac{\chi_{Q_i}}{\sum_k \chi_{Q_k}}$.

By \rf{cc1} we have
$$\mu\left(\bigcup_i 2Q_i\right) \leq \frac{C}{\lambda} \sum_i
\int_{Q_i} |f|\,d\mu \leq \frac{C}{\lambda}\int |f|\,d\mu.$$
So we have to show that
\begin{equation}  \label{eqq}
\mu\Bigl\{x\in \R^d\setminus \bigcup_i 2Q_i:\,|T_\ve f(x)|>\lambda\Bigr\} \leq 
\frac{C}{\lambda}\int |f|\,d\mu.
\end{equation}
Since $\int b_i\,d\mu=0$, $\supp(b_i)\subset R_i$ and
$\|b_i\|_{L^1(\mu)}\leq C\,\int_{Q_i}|f|\,d\mu$, using some standard
estimates we get
$$\int_{\R^d\setminus 2R_i} |T_\ve b_i|\,d\mu \leq C\,\int |b_i|\,d\mu \leq
C\,\int_{Q_i} |f|\,d\mu.$$
Let us see that
\begin{equation}  \label{27b}
\int_{2R_i\setminus 2Q_i} |T_\ve b_i|\,d\mu \leq
C\,\int_{Q_i} |f|\,d\mu
\end{equation}
too. On the one hand, by \rf{cc6} and using the $L^2(\mu)$ boundedness of $T$ 
and that $R_i$ is $(6,6^{n+1})$-doubling we get
\begin{eqnarray*}
\int_{2R_i} |T_\ve \vphi_i| \,d\mu & \leq &
\left(\int_{2R_i} |T_\ve \vphi_i|^2 \,d\mu\right)^{1/2} \, \mu(2R_i)^{1/2} \\
& \leq & C
\left(\int |\vphi_i|^2 \,d\mu\right)^{1/2} \, \mu(R_i)^{1/2} \\
& \leq & C\,\int_{Q_i}|f|\,d\mu.
\end{eqnarray*}
On the other hand, since $\supp(w_i f)\subset Q_i$, if $x\in 2R_i\setminus
2Q_i$, then $|T_\ve f(x)| \leq C\,\int_{Q_i} |f|\,d\mu/|x-x_{Q_i}|^n$,
and so
$$\int_{2R_i\setminus 2Q_i} |T_\ve (w_i\,f)| \,d\mu \leq C\, \int_{
2R_i\setminus 2Q_i} \frac{1}{|x-x_{Q_i}|^n}\,d\mu(x) \times 
\int_{Q_i} |f|\,d\mu,$$
By Lemma \ref{basic},
the first integral on the right hand side is bounded 
by some constant independent of $Q_i$ and $R_i$, since there are no 
$(6,6^{n+1})$-doubling cubes of the form $6^kQ_i$ between $6Q_i$ and $R_i$.
Therefore, \rf{27b} holds.

Then we have
\begin{eqnarray*}
\int_{\R^d\setminus\bigcup_k 2Q_k} |T_\ve b|\,d\mu & \leq &
\sum_i \int_{\R^d\setminus\bigcup_k 2Q_k} |T_\ve b_i|\,d\mu  \\
& \leq & C\,\sum_i \int_{Q_i} |f|\,d\mu \, \leq \, C\,\int |f|\,d\mu.
\end{eqnarray*}
Therefore,
\begin{equation} \label{bb}
\mu\Bigl\{x\in \R^d\setminus \bigcup_i 2Q_i:\,|T_\ve b(x)|>\lambda\Bigr\} \leq 
\frac{C}{\lambda}\int |f|\,d\mu.
\end{equation}
The corresponding integral for the function $g$ is easier to estimate. Taking
into account that $|g|\leq C\,\lambda$, we get
\begin{equation}  \label{gg}
\mu\Bigl\{x\in \R^d\setminus \bigcup_i 2Q_i:\,|T_\ve g(x)|>\lambda\Bigr\} \leq 
\frac{C}{\lambda^2}\int |g|^2\,d\mu \leq 
\frac{C}{\lambda}\int |g|\,d\mu.
\end{equation}
Also, we have
\begin{eqnarray*}
\int |g|\,d\mu & \leq & \int_{\R^d\setminus \bigcup_i Q_i}|f|\,d\mu +
\sum_i \int |\vphi_i|\,d\mu \\
& \leq & \int|f|\,d\mu + \sum_i \int_{Q_i} |f|\,d\mu \, \leq \,
 C\,\int|f|\,d\mu.
\end{eqnarray*} 
Now, by \rf{bb} and \rf{gg} we get \rf{eqq}.
\end{proof}


\section{Proof of Lemma \ref{lem}} \label{secfinal}

{\bf (a)} Taking into account Remark \ref{rem2},
for $\mu$-almost all $x\in\R^d$ such that $|f(x)|>\lambda$, there exists some
cube $Q_x$ satisfying
\begin{equation}  \label{xdr2}
\frac{1}{\mu(2Q_x)} \int_{Q_x} |f| \,d\mu > \frac{\lambda}{2^{d+1}}
\end{equation}
and such that if $Q_x'$ is centered at $x$ with $l(Q_x')>2l(Q_x)$, then
$$\frac{1}{\mu(2Q_x')} \int_{Q_x'} |f| \,d\mu \leq
\frac{\lambda}{2^{d+1}}.$$
Now we can apply Besicovich's covering theorem (see Remark \ref{detall} below) 
to get an almost disjoint subfamily of cubes
$\{Q_i\}_i \subset \{Q_x\}_x$ satisfying \rf{cc1}, \rf{cc2} and \rf{cc3}.

\vspace{3mm}
{\bf(b)} Assume first that the family of cubes $\{Q_i\}_i$ is finite.
Then we may suppose that this family of cubes is ordered in such a way
that the sizes of the cubes $R_i$ are non decreasing (i.e. $l(R_{i+1})\geq
l(R_i)$).
The functions $\vphi_i$ that we will construct will be of the form $\vphi_i
=\alpha_i\,\chi_{A_i}$, with $\alpha_i\in\R$ and $A_i\subset R_i$.
We set $A_1=R_1$ and
$\vphi_1 = \alpha_1\,\chi_{R_1},$
where the constant $\alpha_1$ is chosen so that $\int_{Q_1}f\,w_1\,d\mu=\int
\vphi_1\,d\mu$.

Suppose that $\vphi_1,\ldots,\vphi_{k-1}$ have been constructed,
satisfy \rf{cc4} and $$\sum_{i=1}^{k-1} |\vphi_i|\leq
B\,\lambda,$$  where $B$ is some constant which will be fixed below.

Let $R_{s_1},\ldots,R_{s_m}$ be the subfamily of
$R_1,\ldots,R_{k-1}$ such that $R_{s_j}\cap R_k \neq \varnothing$.
As $l(R_{s_j}) \leq l(R_k)$ (because of the non decreasing sizes of $R_i$),
we have $R_{s_j} \subset 3R_k$. Taking into account that for $i=1,\ldots,k-1$
$$\int |\vphi_i|\,d\mu \leq \int_{Q_i} |f|\,d\mu$$
by \rf{cc4}, and using that $R_k$ is $(6,6^{n+1})$-doubling and \rf{cc2}, we
get
\begin{eqnarray*}
\sum_j \int |\vphi_{s_j}|\,d\mu & \leq & \sum_j
\int_{Q_{s_j}}|f|\,d\mu\\
& \leq & C \int_{3R_k} |f|\,d\mu \, \leq \, C 
\lambda \mu(6R_k)\, \leq \,C_2\lambda\,\mu(R_k).
\end{eqnarray*}
Therefore,
$$\mu\left\{{\textstyle \sum_j} |\vphi_{s_j}| > 2C_2\lambda\right\}\leq
\frac{\mu(R_k)}{2}.$$ So we set
$$A_k = R_k\cap\left\{{\textstyle \sum_j} |\vphi_{s_j}| \leq
2C_2\lambda\right\},$$
and then $\mu(A_k)\geq \mu(R_k)/2.$

The constant $\alpha_k$ is chosen so that for $\vphi_k = \alpha_k\, \chi_{A_k}$
we have $\int\vphi_k\,d\mu = \int_{Q_k} f\,w_k\,
d\mu$.  Then we obtain
$$  
|\alpha_k|  \leq  \frac{1}{\mu(A_k)}\int_{Q_k}  |f|\,d\mu
\leq \frac{2}{\mu(R_k)}\int_{\frac{1}{2}R_k} |f|\,d\mu  \leq C_3\lambda
$$
(this calculation also applies to $k=1$).
Thus,
$$|\vphi_k|+\sum_{j} |\vphi_{s_j}| \leq (2C_2+C_3)\,\lambda.$$
If we choose $B=2C_2+C_3$, \rf{cc5} follows.

Now it is easy to check that \rf{cc6} also holds. Indeed we have
$$\|\vphi_i\|_{L^\infty(\mu)}\, \mu(R_i) 
 \leq C\,|\alpha_i|\,\mu(A_i) 
 =  C\,\left|\int_{Q_i}f\,w_i\,d\mu\right| 
\, \leq \, C\, \int_{Q_i}|f|\,d\mu.
$$

\vspace{3mm}
Suppose now that the collection of cubes $\{Q_i\}_i$ is not finite.
For each fixed $N$ we consider the family of cubes $\{Q_i\}_{1\leq i \leq N}$.
Then, as above, we construct functions $\vphi_1^N,\ldots,\vphi_N^N$ with
$\supp(\vphi_i^N)\subset R_i$ satisfying
$$\int \vphi_i^N \,d\mu = \int_{Q_i} f\,w_i\,d\mu,$$
$$\sum_{i=1}^N |\vphi_i^N| \leq B\,\lambda
$$
and
$$\|\vphi_i^N\|_{L^\infty(\mu)}\, \mu(R_i) \leq C\, \int_{Q_i}|f|\,d\mu.$$
Notice that the sign of $\vphi_i^N$ equals the sign of $\int f\,w_i\,d\mu$ and
so it does not depend on $N$.

Then there is a subsequence
$\{\vphi_1^k\}_{k\in I_1}$ which is convergent in the weak $\ast$ topology of
$L^\infty(\mu)$ to some function
$\vphi_1\in L^\infty(\mu)$. Now we can consider a subsequence
$\{\vphi_2^k\}_{k\in I_2}$ with $I_2\subset I_1$ which
is also convergent in the weak $\ast$ topology of $L^\infty(\mu)$ to some 
function $\vphi_2\in L^\infty(\mu)$.
In general, for each $j$ we consider a subsequence
$\{\vphi_j^k\}_{k\in I_j}$ with $I_j\subset I_{j-1}$ that converges
in the weak $\ast$ topology of $L^\infty(\mu)$ to some function
$\vphi_j\in L^\infty(\mu)$. It is easily checked that the functions
$\vphi_j$ satisfy the required properties.
\qed

\vspace{4mm}
\begin{rem} \label{detall}
Recall that Besicovich's covering theorem asserts that if $\Omega\subset \R^d$
is a {\em bounded} set and for each $x\in\Omega$ there is a cube $Q_x$ centered
at $x$, then there exists a family of cubes $\{Q_{x_i}\}_i$ with finite
overlap covering $\Omega$.

In (a) of the preceeding proof we have applied Besicovich's covering theorem
to $\Omega=\{x:\,|f(x)|>\lambda\}$. However this set may be unbounded,
and the boundedness property is a necessary assumption in Besicovich's theorem
(example: take $\Omega=[0,+\infty)\subset\R$ and consider $Q_x=[0,2x]$ for all
$x\in \Omega$).

We can solve this problem using different arguments. One possibility is to
consider for each $r>0$ the set $\Omega_r = \{x:\,|x|\leq r,\,|f(x)|>\lambda\}$
and to apply Besicovich's covering theorem to $\Omega_r$.
With the same arguments as above, we can decompose $f=g+b$, with $|g|\leq
\lambda$ only on $\Omega_r$ and $b$ as above. Then the proof of Theorem \ref{teo}
can be modified to show that for any fixed constants $\lambda,R>0$ one has
$$\mu\{x\in B(0,R):|T_\ve f(x)|>\lambda\} \leq C\,
\frac{\|f\|_{L^1(\mu)}}{\lambda}.$$

However we prefer the following solution.
We are interested in showing that the Calder\'on-Zygmund decomposition
of Lemma \ref{lem} works also without assuming $\Omega=\{x:\,|f(x)|>\lambda\}$ 
bounded. Let us sketch the argument. Consider a cube $Q_0$ centered at $0$ big
enough so that $$2^{d+1}\,\|f\|_{L^1(\mu)} / \mu(Q_0) <\lambda.$$ 
So for any cube $Q$ containing $Q_0$ we will have
\begin{equation} \label{xdr1}
2^{d+1}\,\|f\|_{L^1(\mu)}/ \mu(Q) <\lambda.
\end{equation}
For $m\geq0$ we set $Q_m:= \left(\frac{5}{4}\right)^m\, Q_0$. For each $m$ we
can apply Besicovich's covering theorem to the annulus $Q_m\setminus Q_{m-1}$
(we take $Q_{-1}:=\varnothing$), with cubes $Q_x$ centered at 
$x\in \supp(\mu)\cap(Q_m\setminus Q_{m-1})$ as in (a) of the proof above,
satisfying \rf{xdr2}.

In this argument we have to be careful with the overlapping among the
cubes belonging to coverings of different annuli. Indeed, there exist some
fixed constants $N$ and $N'$ such that if $m\geq N'$, 
for $x\in \supp(\mu)\cap(Q_m\setminus Q_{m-1})$ we have
\begin{equation} \label{xdr3}
Q_x\subset Q_{m+N}\setminus Q_{m-N}.
\end{equation}
Otherwise, it easily seen that $\ell(Q_x)>\frac{3}{4}\ell(Q_m)$,
choosing $N$ big enough. It follows that $Q_0\subset 2Q_x$ since 
$\ell(Q_0)\ll\ell(Q_m)$ for $N'$ big enough too. This cannot happen because then 
$2Q_x$ satisfies \rf{xdr1}, which contradicts \rf{xdr2}.

Because of \rf{xdr3}, the covering made up of squares belonging to the Besicovich
coverings of different annuli $Q_m\setminus Q_{m-1}$, $m\geq0$, 
will have finite overlap.

Notice that in this argument, it is essential the fact that in \rf{xdr2} we are
not dividing by $\mu(Q_x)$, but by $\mu(2Q_x)$.
\end{rem}


\end{document}